\documentclass[12pt]{amsart}
\usepackage{leftidx}
\usepackage{mathrsfs}
\usepackage{amssymb}
\usepackage{mathtools}
\usepackage{verbatim}
\usepackage{lipsum}
\usepackage{fancyhdr}

\newcommand{\Tr}{\mathrm{Tr}}

\newcommand{\sgn}{\mathrm{sgn}}

\newcommand{\Mp}{\mathrm{Mp}}

\newcommand{\SL}{\mathrm{SL}}

\newtheorem{theorem}{Theorem}
\newtheorem{lemma}{Lemma}
\newtheorem{proposition}{Proposition}
\newtheorem{corollary}{Corollary}

\newtheorem*{example}{Example}

\numberwithin{equation}{section}
\numberwithin{theorem}{section}
\numberwithin{definition}{section}
\numberwithin{lemma}{section}
\numberwithin{proposition}{section}
\numberwithin{corollary}{section}

\newcommand\blfootnote[1]{%
	\begingroup
	\renewcommand\thefootnote{}\footnote{#1}%
	\addtocounter{footnote}{-1}%
	\endgroup
}

\allowdisplaybreaks

\begin{document}
\title{On linear relations for L-values over real quadratic fields}
\pagestyle{fancy}
\rhead[Ren-He Su]{}
\lhead[]{On linear relations for L-values over real quadratic fields}
\author{Ren-He Su}
\address{Department of mathematics, Kyoto University, Kitashirakawa, Kyoto, 606-8502, Japan}
\email{ru-su@math.kyoto-u.ac.jp}
\maketitle

\begin{abstract}
In this paper, we give a method to construct a classical modular form from a Hilbert modular form.
By applying this method, we can get linear formulas which relate the Fourier coefficients of the Hilbert and classical modular forms.
The paper focuses on the Hilbert modular forms over real quadratic fields.
We will state a construction of relations between the special values of L-functions, especially at $0,$ and arithmetic functions.
We will also give a relation between the sum of squares functions with underlying fields $\mathbb{Q}(\sqrt{D})$ and $\mathbb{Q}$.
\end{abstract}

\blfootnote{Acknowledgement. The author would like to show his sincere gratitude to Prof. Ikeda for his insight for the main idea of this paper.}

\section{Introduction}\label{sec:introduction}
This paper aims to give a way to relate the special values of Hecke L-functions over a real quadratic field to arithmetic functions on $\mathbb{Z}$ by restricting the domain of a Hilbert modular form of half-integral weight to construct a classic modular form of odd integral weight.
In \cite{Shinta:76}, Shintani developed a formula for evaluating the special values of Hecke L-functions over totally real number fields at non-positive integers.
Based on the work by Shintani, Okazaki later gave an effective method to evaluate the special values of Hecke L-functions over real quadratic fields by showing how to compute the quadratic characters at ideals explicitly.
However, it seems one can not very easily deduce any linear relation for the L-values from their results. 
\par
A famous formula on relating the L-values to arithmetic functions is Kronecker-Hurwitz class number relation.
For $N\geq0,$ the Hurwitz class number $H(N)$ is defined as follows.
Put $H(0)=-1/12$.
If $N>0$, $H(N)$ is the number of $\SL_2(Z)$-equivalent classes of of not necessarily primitive positive definite binary quadratic forms of discriminant $N$ weighted by the reciprocal of the order of whose stabilizer in $\SL_2(Z)$.
The Hurwitz class number can be written in L-values and arithmetic functions as
\begin{equation}\label{eq:def_of_cohen_h}
H(N)=\begin{cases}
\zeta(-1)&\mbox{if }N=0,\\
L(0,\chi_D)\sum_{d\,|\,f}\mu(d)\chi_D(d)\sigma_1(f/d)&\mbox{if }N>0\mbox{ and }-N=f^2D,\\
0&\mbox{otherwise,}
\end{cases}
\end{equation}
in which $\sigma_1$ is the sum of positive divisors function of degree $1,$ $\chi_D$ is the quadratic character with respect to $D$ and in the second case $f\in\mathbb{N}$ and $D$ is a fundamental discriminant.
Kronecker observed the formula
\[
2\sigma_1(N)=\sum_{s\in\mathbb{Z}}H(4N-s^2)+2\lambda_1(N)
\]
for $N\geq1,$ where
\[
\lambda_1(N)=\sum_{d\,|\,N}\min\{d,N/d\}.
\]
Also, Eichler \cite{Eich:55} showed
\[
\frac{1}{3}\sigma_1(N)=\sum_{s\in\mathbb{Z}}H(N-s^2)+\lambda_1(N)
\]
for odd $N\geq1$.

\par

On 1975, Cohen \cite{HeCo:75} generalized the Hurwitz class numbers to a number $H(r,N),$ which is a special case of equation (\ref{eq:def_of_H}) with underlying field $\mathbb{Q},$ $\chi'=1$ and $\kappa=r$.
We have $H(1,N)=H(N)$.
Similar to the case $r=1,$ the generalized Hurwitz class numbers are closely related to the values of Dirichlet L-functions at $1-2r$.
If we put
\[
\mathscr{H}_{r+1/2}(z)=\sum_{N=0}^\infty H(r,N)q^{2\pi\sqrt{-1}Nz}
\]
for $r\geq2$ and $\Im(z)>0,$ then $\mathscr{H}_{r+1/2}$ is a modular form of weight $r+1/2$.
It is known that $\mathscr{H}_{r+1/2}$ is a Hecke eigenform at all odd primes.
Using the fact that $\mathscr{H}_{r+1/2}$ is a modular form, Cohen derived various formulas such as
\begin{equation}\label{eq:ex1_from_Cohen}
H(2,N)=-\frac{1}{5}\sum_{s\in\mathbb{Z}}\sigma_1\left(\frac{N-s^2}{4}\right)-\begin{cases}
N/10&\mbox{ if }N\mbox{ is a square},\\
0&\mbox{ otherwise}
\end{cases}
\end{equation}
and
\begin{equation}\label{eq:ex2_from_Cohen}
H(4,N)=\sum_{s\in\mathbb{Z}}\sigma_3\left(\frac{N-s^2}{4}\right)
\end{equation}
where we set $\sigma_r(0)=\zeta(-r)/2$.
Also, Cohen showed that if $D\equiv0, 1\mod{4}$ is an integer such that $(-1)^{r+1}D=|D|,$ then for $r\geq2$ we have
\[
\sum_{N=0}^\infty\left(
\sum_{s\in\mathbb{Z}}H\left(
r,\frac{N-s^2}{|D|}
\right)
\right)e^{2\pi\sqrt{-1}Nz}\in M_{r+1}(\Gamma_0(D),\chi_D)
\]
and
\[
\sum_{N=0}^\infty\left(
\sum_{s\in\mathbb{Z},\mathrm{odd}}H\left(
r,\frac{4N-s^2}{|D|}
\right)
\right)e^{2\pi\sqrt{-1}Nz}\in M_{r+1}(\Gamma_0(4D),\chi_D)
\]
where $M_k(\Gamma',\chi)$ denotes the space of modular forms of weight $k$ and the character $\chi$ defined on congruence subgroup $\Gamma'$.
These can be deemed as generalized Kronecker-Hurwitz class number relations.

\par

The modular forms of half integral weight constructed by Cohen was later generalized to the case of Hilbert modular forms by the author in \cite{Ren:14}.
We will in particular employ the generalized version over quadratic field.
Given a real quadratic number field $F$ with discriminant $D>0,$ ring of integers $\mathfrak{o}_F$ and different $\mathfrak{d}_F$ over $\mathbb{Q},$ we put
\[
\gamma\in\Gamma=\left\{
\begin{pmatrix}
a&b\\c&d
\end{pmatrix}\in\SL_2(F)\,\bigg|\,a,d\in\mathfrak{o}_F, b\in\mathfrak{d}_F^{-1}, c\in4\mathfrak{d}_F
\right\}.
\]
If $x\in F,$ we say $x\equiv\square\mod4$ if there exists $\lambda\in\mathfrak{o}_F$ such that $x-\lambda^2\in4\mathfrak{o}_F$ and $x\succ0$ if $x$ is totally positive.
For a character $\chi'$ of the class group of $F$ and integer $\kappa\geq1,$ set
\[
G_{\kappa+1/2,\chi'}
=L_F(1-2\kappa,\overline{\chi'}^2)
+\sum_{
	\begin{smallmatrix}
	(-1)^\kappa\xi\equiv\square\mod4\\
	\xi\succ0
	\end{smallmatrix}}
\mathcal{H}_{\kappa}(\xi,\chi')q^\xi
\]
where $=L_F(s,\overline{\chi'}^2)$ is the L-function over $F$ with character $\overline{\chi'}^2$ and $\mathcal{H}_{\kappa}(\xi,\chi')$ is defined as equation (\ref{eq:def_of_H}), then $G_{\kappa+1/2,\chi'}$ is a Hilbert modular form of parallel weight $\kappa+1/2$ over $F,$ which will be defined explicitly in Section \ref{sec:Results}.
The $q$-coefficients $\mathcal{H}_{\kappa}(\xi,\chi')$ contain information of the Hecke L-values.
Put
\[
\omega=\begin{cases}
(1+\sqrt{D})/2&\mbox{if }D\equiv1\mod{4},\\
\sqrt{D}/2&\mbox{if }D\equiv0\mod{4}.\\
\end{cases}
\] 
We have $\mathfrak{o}_F=\mathbb{Z}+\mathbb{Z}\omega$.
If $u=\alpha+\beta\omega$ is a unit such that $N_{F/\mathbb{Q}}(u)=-1$ and $\alpha, \beta>0,$
then we will show that (Theorem \ref{thm:consts}) for any $n\geq1,$
\[
\sum_{
	\begin{smallmatrix}
	(a,b)\in\mathbb{Z}^2\\
	a\beta-b\alpha=n
	\end{smallmatrix}	
}\mathcal{H}_{\kappa}(a+b\omega,\chi')=\frac{2L_F(1-2\kappa,\overline{\chi'}^2)}{L(-2\kappa,\chi_{-4})}(\sigma_{2\kappa,\chi_{-4}}(n)+(-1)^\kappa\sigma'_{2\kappa,\chi_{-4}}(n))+c(n),
\]
where
\[
\sigma_{2\kappa,\chi_{-4}}(n)=\sum_{d\,|\,n}d^{2\kappa}\chi_{-4}(d),
\]
\[
\sigma'_{2\kappa,\chi_{-4}}(n)=\sum_{d\,|\,n}d^{2\kappa}\chi_{-4}(n/d)
\]
and $c(n)$ is the $q$-coefficient of some cusp form in $S_{2\kappa+1}(\Gamma_0(4),\chi_{-4})$.
In particular, if $\kappa=1,$ then $S_{2\kappa+1}(\Gamma_0(4),\chi_{-4})=0$ thus $c(n)=0$ for all $n$ (Corollary \ref{cor:kappa_is_1}).
\par
We will state the results in Section \ref{sec:Results}, which contains two main theorems.
The first main theorem is the core concept of this paper.
It relates the Hilbert modular forms to modular forms over $\mathbb{Q}$.
The second main theorem, which is what we have just stated above, can be seen as a corollary of the first.
Section \ref{sec:proof_1} and Section \ref{sec:proof_2} are devoted the proofs of the two main theorems.

\section{Results}\label{sec:Results}

Let $F$ be a totally real number field with degree $m$.
The different over $\mathbb{Q},$ ring of integers and real embeddings are denoted by $\mathfrak{d}_F,$ $\mathfrak{o}_F$ and $\iota_i$'s respectively.
Every $x\in F$ will be deemed as a real $n$-tuple $(\iota_1(x),\iota_2(x),\dots,\iota_n(x))\in\mathbb{R}^n$.
As usual, the upper half-plane $\{z=x+\sqrt{-1}y\,|\,x,y\in\mathbb{R},y>0\}$ is denoted by $\mathfrak{h}$.
The standard theta series $\theta_F$ of weight $1/2$ with respect to $F$ is given by
\[
\theta_F(z)=\sum_{\xi\in\mathfrak{o}_F}q^{\xi^2}
\]
where, as usual, for $x\in F$ and $z=(z_1,z_2,\dots,z_m)\in\mathfrak{h}^m,$ we let
\[
q^x=e^{2\pi\sqrt{-1}\sum_{i}z_i\iota_i(x)}.
\]
Then for
\[
\gamma\in\Gamma=\left\{
\begin{pmatrix}
a&b\\c&d
\end{pmatrix}\in\SL_2(F)\,\bigg|\,a,d\in\mathfrak{o}_F, b\in\mathfrak{d}_F^{-1}, c\in4\mathfrak{d}_F
\right\},
\]
and $z\in\mathfrak{h}^m,$ we define the factor of automorphy $j(\gamma,z)$ by
\[
j(\gamma,z)=\frac{\theta_F(\gamma z)}{\theta_F(z)}.
\]
For $\kappa\geq0,$ a Hilbert modular form over $F$ with respect to the congruence subgroup $\Gamma$ and weight $2\kappa+1$ is one corresponding to the factor of automorphy $j(\gamma,z)^{2\kappa+1}$.
The spaces consisting of such Hilbert modular forms with respect to congruence subgroup $\Gamma$ and weight $2\kappa+1$ is denoted by $M_{\kappa+1/2}(\Gamma)$.
We let $S_{\kappa+1/2}(\Gamma)$ be the subspace of $M_{\kappa+1/2}(\Gamma)$ consisting of cusp forms.

\par

Now assume that the different $\mathfrak{d}_F$ is a principal ideal generated by some totally positive element $\boldsymbol{\delta}$.
For $f\in M_{\kappa+1/2}(\Gamma),$ we put
\begin{equation}\label{eq:def_of_Rf}
\mathcal{R}f(z)=f\left(\left(\frac{z}{\iota_1(\boldsymbol{\delta})},\frac{z}{\iota_2(\boldsymbol{\delta})},\dots,\frac{z}{\iota_m(\boldsymbol{\delta})}\right)\right)
\end{equation}
for $z\in\mathfrak{h}$.

\begin{theorem}\label{thm:main}
For $f\in M_{\kappa+1/2}(\Gamma),$ we have
\[
\mathcal{R}f(z)\in M_{m(\kappa+1/2)}(\Gamma_0(4))
\]
where $M_{m(\kappa+1/2)}(\Gamma_0(4))$ is the space of modular forms of weight $m(\kappa+1/2)$ over the underlying field $\mathbb{Q}$.
The factor of automorphy is defined as above for $\mathbb{Q}$
\end{theorem}

We want to focus on real quadratic fields.
Let $F=\mathbb{Q}(\sqrt{D})$ with $D>0$.
We put
\begin{equation}\label{eq:def_of_omega}
\omega=\begin{cases}
(1+\sqrt{D})/2&\mbox{if }D\equiv1\mod{4},\\
\sqrt{D}/2&\mbox{if }D\equiv0\mod{4}.\\
\end{cases}
\end{equation}
Then $\mathfrak{o}_F=\mathbb{Z}+w\mathbb{Z}$.
The different $\mathfrak{d}_F$ is the principal ideal $(\sqrt{D})$ generated by $\sqrt{D}$.
We assume that there exists a unit $u=\alpha+\beta\omega\in\mathfrak{o}_F$ with norm $-1$ such that $\alpha, \beta>0$.
Then $u\sqrt{D}$ is a totally positive integer which generates $\mathfrak{d}_F$.
Putting $\boldsymbol{\delta}=u\sqrt{D},$ we can exactly get the Fourier coefficients of $\mathcal{R}f$ by a straightforward calculation.
We write these results as a corollary of Theorem \ref{thm:main}.

\begin{corollary}\label{cor:quad}
Let $F=Q(\sqrt{D})$ be a real quadratic field with some unit $u=\alpha+\beta\omega$ such that $N_{F/\mathbb{Q}}(u)=-1$ and $\alpha, \beta>0$.
Put $\boldsymbol{\delta}=u\sqrt{D}$.
If we define $\mathcal{R}f$ by (\ref{eq:def_of_Rf}), then it is in $M_{2\kappa+1}(\Gamma_0(4),\chi_{-4}),$ the space of classical modular forms with weight $2\kappa+1$ and the Dirichlet character $\chi_{-4}=\left(\frac{-4}{\cdot}\right)$ associated to $-4$.
Furthermore, if $f=\sum_{\xi\in\mathfrak{o}_F}c(\xi)q^\xi,$ the $q$-expansion for $\mathcal{R}f$ is given by
\[
\mathcal{R}f(z)=\sum_{n=0}^\infty\left(
\sum_{
	\begin{smallmatrix}
	(a,b)\in\mathbb{Z}^2\\
	a\beta-b\alpha=n
	\end{smallmatrix}	
}c(a+b\omega)
\right)q^n.
\]
Here for $z\in\mathfrak{h}$ we put $q^n=e^{2\pi\sqrt{-1}nz}$.
\end{corollary}

By Riemann-Roch theorem, the dimension of $M_{2\kappa+1}(\Gamma_0(4),\chi_{-4})$ and $S_{2\kappa+1}(\Gamma_0(4),\chi_{-4})$ are given by follows:
\begin{align*}
\dim M_{2\kappa+1}(\Gamma_0(4),\chi_{-4})&=0&\quad\mbox{for }\kappa<0;\\
\dim M_{2\kappa+1}(\Gamma_0(4),\chi_{-4})&=1+\kappa&\quad\mbox{for }\kappa\geq0.\\
\dim S_{2\kappa+1}(\Gamma_0(4),\chi_{-4})&=0&\quad\mbox{for }\kappa<2;\\
\dim S_{2\kappa+1}(\Gamma_0(4),\chi_{-4})&=\kappa-1&\quad\mbox{for }\kappa\geq2.
\end{align*}
Using the knowledges about both the Hilbert modular forms and the classical modular forms, one can generate many linear relations between the coefficients of the Hilbert modular forms of half-integral weight and the classical modular forms with odd weight and character $\chi_{-4}$.
For example, let $\theta_F$ be as above and $\theta_\mathbb{Q}=\theta=1+2q+2q^4+2q^9+2q^{16}+\cdots$ be the standard theta function over underlying field $\mathbb{Q}$.
Since $\theta^2\in M_1(\Gamma_0(4),\chi_{-4}),$ which has dimension one, we have $M_1(\Gamma_0(4),\chi_{-4})=\mathbb{C}\cdot\theta^2$.
By comparing the constant term, we get
\begin{equation}\label{sum_of_sqs}
\mathcal{R}\theta_F=\theta^2.
\end{equation}
In particular, if we put
\[
r_{F,k}(\xi)=\#\left\{
(\xi_1,\xi_2,\dots,\xi_k)\in\mathfrak{o}_F^k\ : \xi_1^2+\xi_2^2+\cdots+\xi_k^2=\xi
\right\}
\]
and
\[
r_k(n)=\#\left\{
(n_1,n_2,\dots,n_k)\in\mathbb{Z}^k\ : n_1^2+n_2^2+\cdots+n_k^2=n
\right\},
\]
by equation (\ref{sum_of_sqs}), we get:

\begin{proposition}\label{prop:sum_of_sqs}
Let $\kappa\geq1$ be a positive integer, then
\[
\sum_{\begin{smallmatrix}
	(a,b)\in\mathbb{Z}^2\\
	a\beta-b\alpha=n
	\end{smallmatrix}	
	}r_{F,\kappa}(a+b\omega)=r_{2\kappa}(n).
\]
\end{proposition}

For example, if we take $D=5$ and $u=(1+\sqrt{5})/2,$ the we have
\[
\sum_{m\in\mathbb{Z}}r_{Q(\sqrt{5}),\kappa}\left(n+m\frac{1+\sqrt{5}}{2}\right)=r_{2\kappa}(n),
\]
which is easy to prove.
But for general $D$ and $u,$ it seems that Proposition \ref{prop:sum_of_sqs} is not very obvious.

\par

As mentioned in Section \ref{sec:introduction}, on 1975, Cohen \cite{HeCo:75} constructed a certain modular form $\mathscr{H}_{r+1/2}$ of half-integral weight for the field $\mathbb{Q}$ with respect to $\Gamma_0(4)$.
It is also known as Cohen-Eisenstein series in the Kohnen plus space, which is the subspace of $M_{\kappa+1/2}(\Gamma_0(4))$ which consists of those whose $n$-th Fourier coefficient does not vanish only if $(-1)^\kappa n\equiv0 ,1\mod{4}$.
The subspace of cusp forms in the Kohnen plus space of weight $\kappa+1/2$ was proved to be isomorphic to $S_{2\kappa}(\SL_2(\mathbb{Z}))$ by Kohnen \cite{Kohnen:80}.
The Kohnen plus space was later generalized to the case for Hilbert modular forms of half-integral weight by Hiraga and Ikeda in \cite{HiraIke:13} and the Cohen-Eisenstein series in the corresponding plus space were constructed in \cite{Ren:14}.
As in the $\mathbb{Q}$-case, the Cohen-Eisenstein series in the generalized Kohnen plus space contains information of Hecke L-values and several arithmetic functions in its Fourier coefficients over the underlying real number field.
Thus by Corollary \ref{cor:quad}, we can derive many linear relations between the Hecke L-values over real quadratic fields and the arithmetic functions over $\mathbb{Q}$.

\par

We give the Cohen-Eisenstein series constructed in \cite{Ren:14}.
Let $\kappa$ be a positive integer, which is not $1$ if $F=\mathbb{Q}$, and $\chi'$ a character of the class group of $F$.
For $x\in F,$ we write $x\equiv\square\mod{4}$ if there exists $\lambda\in\mathfrak{o}_F$ such that $x-\lambda^2\in4\mathfrak{o}_F$ and write $x\succ0$ if $x$ is totally positive.
Then the Cohen-Eisenstein series of (parallel) weight $\kappa+1/2$ twisted by $\chi'$ is given by
\[
G_{\kappa+1/2,\chi'}
=L_F(1-2\kappa,\overline{\chi'}^2)
+\sum_{
\begin{smallmatrix}
(-1)^\kappa\xi\equiv\square\mod4\\
\xi\succ0
\end{smallmatrix}}
\mathcal{H}_{\kappa}(\xi,\chi')q^\xi
\]
where
\begin{align}\label{eq:def_of_H}
\mathcal{H}_{\kappa}(\xi,\chi')=&\chi'(\mathcal{D}_{(-1)^\kappa\xi})L_F(1-\kappa,\chi_{(-1)^\kappa\xi}\chi')\nonumber\\
\times&\sum_{\mathfrak{a}\,|\,\mathfrak{f}_{(-1)^\kappa\xi}}\mu_F(\mathfrak{a})\chi_{(-1)^\kappa\xi}(\mathfrak{a})\chi'(\mathfrak{a})N_{F/\mathbb{Q}}(\mathfrak{a})^{\kappa-1}\sigma_{F,2\kappa-1,{\chi'}^2}(\mathfrak{f}_{(-1)^\kappa\xi}\mathfrak{a}^{-1}).
\end{align}
Here $\mathcal{D}_x$ and $\chi_x$ are the relative discriminant and the quadratic character, respectively, associated to the quadratic extension $F(\sqrt{x})/F$.
$\mathfrak{f}_x$ is the integral ideal such that $\mathfrak{f}_x^2\mathcal{D}_x=(x)$.
In the sum in equation \ref{eq:def_of_H}, the ideal $\mathfrak{a}$ runs over all integral ones dividing $\mathfrak{f}_{(-1)^\kappa\xi}$.
Finally, $\mu_F$ is the M\"{o}bius function and 
\[
\sigma_{F,k,\chi}(\mathfrak{b})=\sum_{\mathfrak{r}\,|\,\mathfrak{b}}N_{F/Q}(\mathfrak{r})^k\chi(\mathfrak{r}).
\]
This form $G_{\kappa+1/2,\chi'}$ is a Hecke eigenform in $M_{\kappa+1/2}(\Gamma)$ at any odd non-archimedean place $v$ of $F$.

\par

We can apply Corollary \ref{cor:quad} on the Cohen-Eisenstein series for real quadratic fields.
Let $\sigma_{k,\chi_{-4}}(n)=\sigma_{\mathbb{Q},k,\chi_{-4}}(n)$ and
\[
\sigma'_{k,\chi_{-4}}(n)=\sum_{d\,|\,n}d^k\chi_{-4}(n/d).
\]
For $\kappa\geq1,$ two normalized Eisenstein series in $M_{2\kappa+1}(\Gamma_0(4),\chi_{-4})$ are given by
\[
E_{2\kappa+1,\chi_{-4}}(z)=\frac{1}{2}L(-2\kappa,\chi_{-4})+\sum_{n=1}^\infty\sigma_{2\kappa,\chi_{-4}}(n)q^n
\]
and
\[
F_{2\kappa+1,\chi_{-4}}(z)=\sum_{n=1}^\infty\sigma'_{2\kappa,\chi_{-4}}(n)q^n.
\]
They span the subspace of Eisenstein series in $M_{2\kappa+1}(\Gamma_0(4),\chi_{-4}),$ that is, the orthogonal complement of $S_{2\kappa+1}(\Gamma_0(4),\chi_{-4})$.
The coefficients of $E_{2\kappa+1,\chi_{-4}}$ and $F_{2\kappa+1,\chi_{-4}}$ in the linear combination for $\mathcal{R}G_{\kappa+1/2,\chi'}$ can be explicitly calculated as follows.

\begin{theorem}\label{thm:consts}
Let $F$ be a real quadratic field with a unit with norm $-1$.
We have
\[
\mathcal{R}G_{\kappa+1/2,\chi'}=\frac{2L_F(1-2\kappa,\overline{\chi'}^2)}{L(-2\kappa,\chi_{-4})}(E_{2\kappa+1,\chi_{-4}}+(-1)^\kappa F_{2\kappa+1,\chi_{-4}})+Q
\]
where $Q\in S_{2\kappa+1}(\Gamma_0(4),\chi_{-4})$ is a cusp form.
\end{theorem}

In particular, for $\kappa=1,$ since $S_{3}(\Gamma_0(4),\chi_{-4})=0$ and $L(-2,\chi_{-4})=-1/2,$ the next corollary follows.

\begin{corollary}\label{cor:kappa_is_1}
With the notations above, for any $n\geq1,$ we have
\[
\sum_{
	\begin{smallmatrix}
	(a,b)\in\mathbb{Z}^2\\
	a\beta-b\alpha=n
	\end{smallmatrix}	
}\mathcal{H}_{1}(a+b\omega,\chi')=-4L_F(1-2\kappa,\overline{\chi'}^2)(\sigma_{2,\chi_{-4}}(n)-\sigma'_{2,\chi_{-4}}(n)).
\]	
\end{corollary}

\begin{example}
If we put
\[
\mathcal{S}_5(z)=\sum_{n=0}^\infty s(n)q^n
\]
where
\[
s(n)=\frac{1}{4}\sum_{a^2+b^2=n}(a+b\sqrt{-1})^4,
\]
then $\mathcal{S}_5$ is a cusp form in $S_5(\Gamma_0(-4),\chi_{-4})$.
Thus $M_5(\Gamma_0(4),\chi_{-4})$ is spanned by $E_{5,\chi_{-4}},$ $F_{5,\chi_{-4}}$ and $\mathcal{S}_5$ over $\mathbb{C}$.

\par

Set $F=Q(\sqrt{5}).$
We let $u=\omega=(1+\sqrt{5})/2$.
By some calculations, we have
\[
\mathcal{R}(G_{3/2,1})=-\frac{2}{15}(E_{3,\chi_{-4}}-F_{3,\chi_{-4}})
\]
and
\[
\mathcal{R}(G_{5/2,1})=\frac{1}{75}(E_{5,\chi_{-4}}+F_{5,\chi_{-4}})+\frac{1}{25}\mathcal{S}_5.
\]
These yield 

\begin{align*}
L_F(0,\chi_{-2-\omega})&=-\frac{2}{15}\left(\sigma_{2,\chi_{-4}}(2)-\sigma'_{2,\chi_{-4}}(2)\right),\\
2L_F(0,\chi_{-3})+2L_F(0,\chi_{-3+\omega})&=-\frac{2}{15}\left(\sigma_{2,\chi_{-4}}(3)-\sigma'_{2,\chi_{-4}}(3)\right),\\
2L_F(0,\chi_{-4})&=-\frac{2}{15}\left(\sigma_{2,\chi_{-4}}(4)-\sigma'_{2,\chi_{-4}}(4)\right),\\
2L_F(0,\chi_{-3})+2L_F(0,\chi_{-6-\omega})&=-\frac{2}{15}\left(\sigma_{2,\chi_{-4}}(6)-\sigma'_{2,\chi_{-4}}(6)\right),\\
2\zeta_F(-1)&=\frac{1}{75}\left(\sigma_{4,\chi_{-4}}(1)+\sigma'_{4,\chi_{-4}}(1)+3s(1)\right),\\
4\zeta_F(-1)+2L_F(-1,\chi_{5+\omega})&=\frac{1}{75}\left(\sigma_{4,\chi_{-4}}(5)+\sigma'_{4,\chi_{-4}}(5)+3s(5)\right),...
\end{align*}
and so on.
\end{example}

\par

\section{Proof for Theorem \ref{thm:main}}\label{sec:proof_1}

Fix $f\in M_{\kappa+1/2}(\Gamma)$.
For $\gamma\in\Gamma_0(4)\subset\SL_2(\mathbb{Z})$ and $z\in\mathfrak{h},$ by an easy calculation, we have
\[
\mathcal{R}f(\gamma z)=\left(\frac{\mathcal{R}\theta_F(\gamma z)}{\mathcal{R}\theta_F(z)}\right)^{2\kappa+1}\mathcal{R}f(z),
\]
thus we only need to show
\begin{equation}\label{eq:automorphy}
\frac{\mathcal{R}\theta_F(\gamma z)}{\mathcal{R}\theta_F(z)}=\left(\frac{\theta(\gamma z)}{\theta(z)}\right)^m.
\end{equation}
It is well-known that $\Gamma_0(4)$ is generated by the three elements $\begin{pmatrix}
1&1\\0&1
\end{pmatrix}, \begin{pmatrix}
1&0\\4&1
\end{pmatrix}$ and $\begin{pmatrix}
-1&0\\0&-1
\end{pmatrix}$.
Notice
\[
\begin{pmatrix}
1&0\\4&1
\end{pmatrix}=\begin{pmatrix}
0&-2^{-1}\\2&0
\end{pmatrix}\begin{pmatrix}
1&-1\\0&1
\end{pmatrix}\begin{pmatrix}
0&2^{-1}\\-2&0
\end{pmatrix}
\] and 
\[
\begin{pmatrix}
-1&0\\0&-1
\end{pmatrix}=\begin{pmatrix}
0&-2^{-1}\\2&0
\end{pmatrix}^2.
\]
Also, it is well-known that
\[
\frac{\theta(z+1)}{\theta(z)}=1
\]
and
\[
\frac{\theta(-1/4z)}{\theta(z)}=\sqrt{\frac{2z}{\sqrt{-1}}}.
\]
Therefore, to show equation (\ref{eq:automorphy}), it suffices to prove
\begin{equation}\label{eq:pf1}
\frac{\mathcal{R}\theta_F(z+1)}{\mathcal{R}\theta_F(z)}=1
\end{equation}
and
\begin{equation}\label{eq:pf2}
\frac{\mathcal{R}\theta_F(-1/4z)}{\mathcal{R}\theta_F(z)}=\left(\frac{2z}{\sqrt{-1}}\right)^{m/2}.
\end{equation}
The first condition is trivial, so let us prove the second condition (\ref{eq:pf2}).
\par
We first recall the famous Poisson summation formula, for which one can consult for example \cite{AnWe:74}.
\begin{lemma}
Let $g$ be a rapidly decreasing function on $\mathbb{R}^m$ such that $G(t)=\sum_{\xi\in\mathfrak{o}_F}g(\xi+t)$ for $t\in\mathbb{R}^m$ is also rapidly decreasing on $\mathbb{R}^m,$ then we have
\[
\sum_{\xi\in\mathfrak{o}_F}g(\xi)=N_{F/\mathbb{Q}}(\boldsymbol{\delta})^{-1/2}\sum_{\xi\in\mathfrak{d}_F^{-1}}\hat{g}(\xi)
\]
where
\[
\hat{g}(t)=\int_{\mathbb{R}^m}g(s)e^{2\pi\sqrt{-1}\sum_{i=1}^mt_is_i}ds
\]
for $t=(t_1,t_2,\dots,t_m)\in\mathbb{R}^m$.
\end{lemma}
Now by this lemma, we have
\begin{align*}
 &\mathcal{R}\theta_F(-1/4z)\\
=&\sum_{\xi\in\mathfrak{o}_F}\exp\left(-\frac{\pi\sqrt{-1}}{2z}\Tr(\xi^2/\boldsymbol{\delta})\right)\\
=&N_{F/\mathbb{Q}}(\boldsymbol{\delta})^{-1/2}\sum_{\xi\in\mathfrak{d}_F^{-1}}\int_{\mathbb{R}^m}\exp\left(-\frac{\pi\sqrt{-1}}{2z}\sum_{i=1}^m\frac{s_i^2}{\iota_i(\boldsymbol{\delta})}
+2\pi\sqrt{-1}\sum_{i=1}^m s_i\iota_i(\xi)\right)ds\\
=&N_{F/\mathbb{Q}}(\boldsymbol{\delta})^{-1/2}\sum_{\xi\in\mathfrak{d}_F^{-1}}\exp(2\pi\sqrt{-1}z\Tr(\boldsymbol{\delta}\xi^2))\int_{\mathbb{R}^m}\exp\left(-\frac{\pi\sqrt{-1}}{2z}\sum_{i=1}^m\frac{s_i^2}{\iota_i(\boldsymbol{\delta})}\right)ds\\
=&N_{F/\mathbb{Q}}(\boldsymbol{\delta})^{-1/2}\left(\prod_{i=1}^m\sqrt{\frac{2\iota_i(\boldsymbol{\delta})z}{\sqrt{-1}}}\right)\sum_{\xi\in\mathfrak{o}_F}\exp(2\pi\sqrt{-1}z\Tr(\xi^2/\boldsymbol{\delta}))\\
=&\left(\frac{2z}{\sqrt{-1}}\right)^{m/2}\mathcal{R}\theta_F(z)
\end{align*}
where in the third equation we shifted the integration paths and in the fourth equation we used the fact that $\boldsymbol{\delta}$ is totally positive.
This concludes our proof.

\section{Proof for Theorem \ref{thm:consts}}\label{sec:proof_2}

Let the Fricke involution $\mathcal{W}$ on $M_{2\kappa+1}(\Gamma_0(4),\chi_{-4})$ be defined by
\[
\mathcal{W}h(z)=\sqrt{-1}(2z)^{-2\kappa-1}h(-1/4z).
\]
Since
\[
E_{2\kappa+1,\chi_{-4}}(z)=(-1)^\kappa\frac{2^{2\kappa-1}(2\kappa)!}{\pi^{2\kappa+1}}\sum_{(0,0)\neq(m,n)\in\mathbb{Z}^2}\frac{\chi_{-4}(n)}{(4mz+n)^{2\kappa+1}}
\]
and
\[
F_{2\kappa+1,\chi_{-4}}(z)=\frac{(2\kappa)!}{2(-2\pi\sqrt{-1})^{2\kappa+1}}\sum_{(0,0)\neq(m,n)\in\mathbb{Z}^2}\frac{\chi_{-4}(m)}{(mz+n)^{2\kappa+1}},
\]
which can be get from applying the Poisson summation formula on the right hand sides, one easily deduces that
\[
\mathcal{W}E_{2\kappa+1,\chi_{-4}}=4^\kappa F_{2\kappa+1,\chi_{-4}}.
\]
Since $\mathcal{W}$ is an involution, one also has
\[
\mathcal{W}F_{2\kappa+1,\chi_{-4}}=4^{-\kappa}E_{2\kappa+1,\chi_{-4}}.
\]

\par

Let
\[
\mathcal{R}G_{\kappa+1/2,\chi'}=c_1E_{2\kappa+1,\chi_{-4}}+c_2F_{2\kappa+1,\chi_{-4}}+Q
\]
where $Q$ is a cusp form.
Notice that in this linear combination, the Eisenstein series $E_{2\kappa+1,\chi_{-4}}$ is the only one which has a nonzero constant term in whose $q$-expansion.
Therefore one first see immediately that
\[
c_1=2L_F(1-2\kappa,\overline{\chi'}^2)L(-2\kappa,\chi_{-4})^{-1}.
\]
Since 
\[
\mathcal{W}\mathcal{R}G_{\kappa+1/2,\chi'}=4^\kappa c_1F_{2\kappa+1,\chi_{-4}}+4^{-\kappa}c_2E_{2\kappa+1,\chi_{-4}}+\mathcal{W}Q
\]
and $\mathcal{W}Q$ is still a cusp form, to prove the theorem, it amounts to show that the constant term of $\mathcal{W}\mathcal{R}G_{\kappa+1/2,\chi'}$ is $(-4)^{-\kappa}L_F(1-2\kappa,\overline{\chi'}^2)$.
Put
\[
\mathcal{W}_FG_{\kappa+1/2,\chi'}(z)=\sqrt{-1}(4Dz_1z_2)^{-\kappa-1/2}G_{\kappa+1/2,\chi'}\left(\left(-(4\boldsymbol{\delta}^2z_1)^{-1},-(4\bar{\boldsymbol{\delta}}^2z_2)^{-1}\right)\right)
\]
for $z=(z_1,z_2)\in\mathfrak{h}^2.$
A simple calculation shows that
\[
\mathcal{RW}_FG_{\kappa+1/2,\chi'}=\mathcal{W}\mathcal{R}G_{\kappa+1/2,\chi'}.
\]
Note that $\mathcal{R}$ does not affect the constant term in the $q$-expansion.
Thus we need to show that $\mathcal{W}_FG_{\kappa+1/2,\chi'}$ has constant term $(-4)^{-\kappa}L_F(1-2\kappa,\overline{\chi'}^2),$ that is, $(-4)^{-\kappa}$ times that of $G_{\kappa+1/2,\chi'}$.
\par
In the paper \cite{Ren:14}, the Eisenstein series was originally defined as an automorphic form.
In order to describe the proof, we gives a brief introduction for the automorphic forms.
One can consult \cite{HiraIke:13} for more details.
\par
For any local place $v$ of $F$, let $F_v$ be the local component of $F$ with respect to $v$.
The metaplectic double covering of $\SL_2(F_v)$ is denoted by $\Mp_2(F_v)$.
Any element in $\Mp_2(F_v)$ has the form $[g,\zeta]$ where $g\in\SL_2(F_v)$ and $\zeta\in\{\pm1\}$.
If $v$ is an odd finite place, there exists a canonical splitting contained in $\Mp_2(F_v)$ over the standard maximal subgroup $K_v$ of $\SL_2(F_v),$ which we also denote by $K_v$.
Let $\mathbb{A}_F$ be the adele ring of $F$.
Then the global metaplectic group $\Mp_2(\mathbb{A}_F)$ is the restricted product $\bigotimes'_v\Mp_2(F_v)$ with respect to $\{K_v\}$ divided by $\{\oplus_v\zeta_v\,|\,\zeta_v\in\{\pm1\},\prod_v\zeta_v=1\}$.
\par
It is known that $\SL_2(F)$ can be uniquely embedded into $\Mp_2(\mathbb{A}_F)$.
The image of this embedding is also denoted by $\SL_2(F)$.
There is a one-to-one corresponding between the modular forms in $M_{\kappa+1/2}(\Gamma)$ and certain kind of automorphic forms on $\SL_2(F)\backslash\Mp_2(\mathbb{A}_F)$.
We state briefly how this correspondence works.
\par
Let $\tilde{j}$ be the unique factor of automorphy on $\Mp_2(\mathbb{R})\times\mathfrak{h}$ such that $\tilde{j}([g,\zeta],\tau)^2=c\tau+d$ for $g=\begin{pmatrix}
a&b\\c&d
\end{pmatrix}\in\SL_2(\mathbb{R})$.
There exists a genuine character $\varepsilon_\mathrm{f}$ on $\widetilde{\prod_{v<\infty}\Gamma_v}$ such that
\[
(\theta_F(\gamma z)/\theta_F(z))^{2\kappa+1}=\varepsilon_\mathrm{f}([\gamma,1])\prod_{i=1}^2\tilde{j}([\iota_i(\gamma),1],z_i)^{2\kappa+1}
\]
where $\Gamma_v\subset\SL_2(F_v)$ is the local component of $\Gamma$ with respect to the place $v<\infty$ and $\widetilde{\prod_{v<\infty}\Gamma_v}$ is the inverse image of $\prod_{v<\infty}\Gamma_v$ in $\Mp_2(\mathbb{A}_F)$.
\par
Fix $f\in M_{\kappa+1/2}(\Gamma)$
We can associate $f$ to a automorphic form $\phi_f$ on $\SL_2(F)\backslash\Mp_2(\mathbb{A}_F)$ by the followings.
For any $\tilde{g}\in\Mp_2(\mathbb{A}_F),$ by the strong approximation theorem, there exist $\gamma\in\SL_2(F),$ $\tilde{g}_\mathrm{f}\in\widetilde{\prod_{v<\infty}\Gamma_v}$ and $\tilde{g}_\infty\in\widetilde{\SL_2(\mathbb{R})^2}$ such that $\tilde{g}=\gamma\tilde{g}_{\mathrm{f}}\tilde{g}_\infty$.
Note that $\tilde{g}_\infty=(\tilde{g}_\mathrm{\infty,1},\tilde{g}_\mathrm{\infty,2})$ acts on $\mathbf{i}=(\sqrt{-1},\sqrt{-1})\in\mathfrak{h}^2$ as the usual M\"{o}bius transformation.
The value of $\phi_f$ at $\tilde{g}$ is defined by
\[
\phi_f(\tilde{g})=f(\tilde{g}_\infty(\mathbf{i}))\left(\varepsilon_\mathrm{f}(\tilde{g}_\mathrm{f})\prod_{i=1}^2\tilde{j}(\tilde{g}_{\infty,i},\sqrt{-1})^{2\kappa+1}\right)^{-1}.
\]
This definition does not depend on the choice of the decomposition of $\tilde{g}$.
We put
\[
\mathcal{A}_{\kappa+1/2}(\SL_2(F)\backslash\Mp_2(\mathbb{A}_F);\varepsilon)=\{
\phi_f\,|\,f\in M_{\kappa+1/2}(\Gamma)
\}.
\]
Conversely, for any $\phi\in\mathcal{A}_{\kappa+1/2}(\SL_2(F)\backslash\Mp_2(\mathbb{A}_F);\varepsilon),$ let
\begin{equation}\label{eq:corresp}
f_\phi(z)=\phi(\tilde{g}_\infty)\prod_{i=1}^2\tilde{j}(\tilde{g}_{\infty,i}\sqrt{-1})^{2\kappa+1}, z\in\mathfrak{h}^2,
\end{equation}
where $\tilde{g}_\infty=(\tilde{g}_\mathrm{\infty,1},\tilde{g}_\mathrm{\infty,2})\in\widetilde{\SL_2(\mathbb{R})^2}$ is chosen so that $\tilde{g}_\infty(\mathbf{i})=z$.
Then this gives an inverse for the lifting $f\mapsto\phi_f$.
Thus we get the one-to-one correspondence between $M_{\kappa+1/2}(\Gamma)$ and $\mathcal{A}_{\kappa+1/2}(\SL_2(F)\backslash\Mp_2(\mathbb{A}_F);\varepsilon)$.
\par
Back to our proof for Theorem \ref{thm:consts}.
Let $B$ be the subgroup of $\SL_2(F)$ consisting of upper triangular elements and $f=\prod_{v\geq\infty}f_v$ be the function on $\Mp_2(\mathbb{A}_F)$ defined as equation (8.1) in \cite{Ren:14}.
Then $G_{\kappa+1/2,\chi'}$ corresponds to an automorphic form $cE'$ where $c$ is a constant and
\[
E'(g)=\sum_{\gamma\in B\backslash\SL_2(F)}f(\gamma g)
\]
is an automorphic form on $\SL_2(F)\backslash\Mp_2(\mathbb{A}_F)$.
We let $\mathbf{w}_{-2\boldsymbol{\delta},\mathrm{f}}$ be the finite part of the image of $\begin{pmatrix}0&(2\boldsymbol{\delta})^{-1}\\-2\boldsymbol{\delta}&0\end{pmatrix}$ in $\Mp_2(\mathbb{A}_F)$.
By equation (\ref{eq:corresp}), the automorphic form on $\SL_2(F)\backslash\Mp_2(\mathbb{A}_F)$ corresponding to $\mathcal{W}_FG_{\kappa+1/2,\chi'}$ is $\sqrt{-1}c\rho(\mathbf{w}_{-2\boldsymbol{\delta},\mathrm{f}})E',$ where $\rho$ is the right translation.
Denote the modular forms corresponding to $E'$ and $\rho(\mathbf{w}_{-2\boldsymbol{\delta},\mathrm{f}})E'$ by $E$ and $\rho(\mathbf{w}_{-2\boldsymbol{\delta},\mathrm{f}})E,$ respectively.
It remains to compare the constant terms of $E$ and $\sqrt{-1}\rho(\mathbf{w}_{-2\boldsymbol{\delta},\mathrm{f}})E$.
\par
The finite part of $f$ is $\prod_{v<\infty}f^+_{K,v},$ where $f^+_{K,v}$ is a function on $\Mp_2(F_v)$ defined as Definition 2.4 in \cite{Ren:14}.
For $E,$ the calculations in page 712 (for $\kappa\geq2$) and page 717 (for $\kappa=1$) of \cite{Ren:14} show that the constant term of $E$ is 
\[
\prod_{v<\infty}f^+_{K,v}(1)=N_{F/\mathbb{Q}}(2)^{2\kappa-1/2}N_{F/\mathbb{Q}}(\boldsymbol{\delta})^{\kappa-1/2}\prod_{v<\infty}\alpha_v((-1)^\kappa),
\]
where $\alpha_v$ is the Weil constant with respect to the place $v$, an eighth root of $1$.
On the other hand, via similar calculations, the constant term of $\sqrt{-1}\rho(\mathbf{w}_{-2\boldsymbol{\delta},\mathrm{f}})E$ is $\sqrt{-1}\prod_{v<\infty}f^+_{K,v}(\mathbf{w}_{-2\boldsymbol{\delta},\mathrm{f}})$.
Then by Proposition 4.3 of \cite{Ren:14}, we have
\[
\sqrt{-1}\prod_{v<\infty}f^+_{K,v}(\mathbf{w}_{-2\boldsymbol{\delta},\mathrm{f}})=\sqrt{-1}N_{F/\mathbb{Q}}(2\boldsymbol{\delta})^{\kappa-1/2}\prod_{v<\infty}\alpha_v((-1)^\kappa)\alpha_v((-1)^\kappa2\boldsymbol{\delta}).
\]
For any $x\in F,$ the Weil constants $\alpha_v$ have the property
\[
\prod_{v\geq\infty}\alpha_v(x)=1.
\]
Thus we have
\[
\sqrt{-1}\prod_{v<\infty}f^+_{K,v}(\mathbf{w}_{-2\boldsymbol{\delta},\mathrm{f}})=\sqrt{-1}N_{F/\mathbb{Q}}(2)^{-\kappa}\prod_{v\mid\infty}\overline{\alpha_v((-1)^\kappa2\boldsymbol{\delta})}\cdot\prod_{v<\infty}f^+_{K,v}(1).
\]
The value of $\alpha_v$ for real places $v$ of $F$ is given by
\[
\alpha_v(x)=\exp(\sgn(x)\pi\sqrt{-1}/4), x\in\mathbb{R}.
\]
Since we have taken $\boldsymbol{\delta}\succ0$ and $F$ is a quadratic field, we get
\[
\sqrt{-1}\prod_{v<\infty}f^+_{K,v}(\mathbf{w}_{-2\boldsymbol{\delta},\mathrm{f}})=(-4)^{-\kappa}\prod_{v<\infty}f^+_{K,v}(1).
\]
Thus we have shown that $\mathcal{W}_F$ modifies the constant term of $G_{\kappa+1/2,\chi'}$ by a $(-4)^{-\kappa}$-multiplication .
This concludes our proof.

\end{document}